\documentclass[leqno,12pt]{article}
\usepackage{amsmath, amssymb}
\usepackage[dvips]{graphicx}
\usepackage{amsthm}
\usepackage{color}%\UTF{0092}\UTF{00C7}\UTF{0089}\UTF{00C1}%
\usepackage{ulem}%\UTF{0092}\UTF{00C7}\UTF{0089}\UTF{00C1}%
\theoremstyle{plain} % イタリック体
\newtheorem{theorem}{\indent\sc Theorem}[section] % 見出しはスモールキャップ
\newtheorem{lemma}[theorem]{\indent\sc Lemma}
\newtheorem{corollary}[theorem]{\indent\sc Corollary}
\newtheorem{proposition}[theorem]{\indent\sc Proposition}

\theoremstyle{definition} % ローマン体に変更

\newtheorem{remark}[theorem]{\indent\sc Remark}

\title{Complete parallel mean curvature surfaces in two-dimensional complex space-forms }

\author{ Katsuei  Kenmotsu }

\makeatletter

\@addtoreset{equation}{section}
\makeatother

\date{}

\begin{document}
\maketitle

\begin{abstract}
 The purpose of this article is to determine explicitly  the  complete surfaces with parallel mean curvature vector, both  in the complex projective plane  and the complex hyperbolic plane.
The main results  are as follows: When the curvature of the ambient space is positive, there exists a unique such  surface up to rigid motions of the target space. On the other hand, when the curvature of  the ambient space is  negative, there are `non-trivial' complete parallel mean curvature surfaces generated by Jacobi elliptic functions and they exhaust such surfaces.
\end{abstract}
\section{Introduction}
 Recently, the study of non-zero parallel mean curvature surfaces with  codimension two has progressed   in the various ambient spaces  (see  Manzano, Torralbo, and Van der Veken \cite{manzano} for references).  Such surfaces in complex space-forms are already classified when the Kaehler angle  is constant by  Chen \cite{chen3}, and finally by Hirakawa \cite{hirakawa}. %, and Fetcu \cite{fetcu}. 
If the Kaehler angles  of these surfaces are not constant,  Ogata \cite{ogata} was the first to study such surfaces in  non-flat complex space-forms, and later  in Kenmotsu and Zhou \cite{kenzhou}, Kenmotsu \cite{ken1},  Hirakawa \cite{hirakawa}, Ferreira and Tribuzy \cite{fertri},  and Fetcu \cite{fetcu}, they continued the research.  
 In particular,  Kenmotsu and Zhou \cite{kenzhou} and Kenmotsu \cite{ken1} proved that  the first and second fundamental forms of such a surface in non-flat complex space-forms are determined by a harmonic function and five real constants. We shall note that a parallel mean curvature surface in a flat complex space-form with codimension two  is contained in a totally umbilical real hypersurface of the ambient space as a minimal or  constant mean curvature surface (Chen \cite{chen1}, Yau \cite{yau}), and also if the mean curvature vector   vanishes identically,  there are many articles studying such surfaces. (see, for instance, Eschenburg, Guadalupe, and Tribuzy \cite{escgt}.) 

In this article, by the refinement of our  method developed in Kenmotsu and Ogata \cite{kenoga} and Kenmotsu \cite{ken2}, we prove that  any parallel mean curvature surface in the complex projective plane  $\mathbb{C}P^{2}$ with the Fubini-Study metric  must have  constant Kaehler angle, thus it  is congruent to  the flat torus if the surface is complete.   For parallel mean curvature surfaces in the complex hyperbolic plane  $\mathbb{C}H^{2}$, with its canonical structure of Kaehler surface, the situation is different from  $\mathbb{C}P^{2}$-case.  There is a two-parameter family of those surfaces, say   $x_{t,p}, \ (0 \leq t \leq \pi,  0 < p\; (\neq 1/4) < \infty)$, such that any parallel mean curvature surface of a general type in $\mathbb{C}H^{2}$  is locally congruent to some $x_{t,p}$. Moreover, there exists the limit  of $x_{t,p}$, as $p$ tends to $0$,  that is  Hirakawa surface, and also letting $p \rightarrow 1/4+$ and  $p \rightarrow 1/4-$, there are the limits of $x_{t,p}$  which are the same ones presented in Section 4 of Kenmotsu  \cite{ken1}. %Therefore, the moduli space of parallel mean curvature surfaces with non-constant Kaehler angle is a rectangle in $R^{2}$ that is $[0,\pi] \times [0,1/3]$.

 In Section 2, we refine  Lemma 3.3 of Kenmotsu and Zhou \cite{kenzhou}, because  its statement  is incomplete. In Section 3, the key result of this article is proved by using our previous results  \cite{ken1}, \cite{ken2} that is Lemma 3.4: a constant appeared in the second fundamental form in \cite{ken2} must be zero if the ambient space is non-flat. In Section 4, we will concentrate on the study in $\mathbb{C}H^{2}$-case and prove that any parallel mean curvature surface of a general type in $\mathbb{C}H^{2}$ is  determined by two real constants, one of which is the parameter of the associated family, and the other  is  the modulus of a Jacobi  elliptic function of the first kind.  The moduli space of  the complete parallel mean curvature surfaces in $\mathbb{C}H^{2}$ is presented at the end of this section.  In Section 5 we state some open questions related to this article.

\section{Parallel mean curvature surfaces of a special type}
This section  establishes the notation that will be used in this article, recalls the structure equations of parallel mean curvature surfaces, and defines two families of parallel mean curvature  surfaces as a refinement of Lemma 3.3 of  \cite{kenzhou}.

Let $\overline{M}[4\rho]$ be a complex  two-dimensional complex space-form with  constant holomorphic sectional curvature $4\rho$,  and let
$M$ be an oriented and connected real two-dimensional Riemannian manifold with  Gaussian curvature
$K$ and $x: M\longrightarrow
\overline{M}[4\rho]$ be an isometric immersion, with   Kaehler angle  $\alpha$ such that the  mean
curvature vector  $H$  is nonzero and parallel for the normal connection on the normal bundle
of the immersion.  Since the length of the mean curvature vector is constant,  we write  $|H| =2b >0$. Hereafter, we call such an immersion as a parallel mean curvature surface.

Let $M_{0} = \{p \in M\ |\  x \ \mbox{is neither holomorphic nor anti-holomorphic at} \  p \}$.  $M_{0}$ is an open dense subset of $M$. Because all of the  calculations and formulas on $M_{0}$ presented  in  \cite{ogata} are valid until page 400, according to a remark in  \cite{hirakawa},  there exists a local
field of unitary coframes
$\{w_{1},w_{2}\}$ on $M_{0}$ such that, by restricting it to $x$, the Riemannian metric $ds^2$
 on $M_{0}$ can be written as  $ds^{2}=\phi\bar{\phi}$, where
$\phi=\cos\alpha/2\cdot\omega_{1}+\sin\alpha/2\cdot\bar{\omega}_{2}$ % (see \cite{ken1})
.  Let $a$ and $c$ be the complex-valued functions on $M_{0}$ that determine the second fundamental form of $x$.
 Then  the Kaehler angle $\alpha$ and the complex 1-form $\phi$  satisfy
\begin{equation}
d\alpha = ( a +  b )\phi + (\bar{a} + b )\bar{\phi} , \quad \mbox{and} \quad
d\phi =   (\bar{a} - b)\cot \alpha \, \phi \wedge \bar{\phi}.
\end{equation}
%where $2b = |H| >0$.%!
 Equations  (2.4), (2.5), and (2.6)  in \cite{ogata} are
 \begin{eqnarray}
&& K= -4(|a|^{2} - b^{2})+6\rho\cos^{2}\alpha, \\
&& da\wedge\phi= -\left(2a (\bar{a} - b) \cot\alpha  + \frac{3}{2}\rho\sin \alpha \cos \alpha \right)
\phi\wedge\bar{\phi}, \\
&& dc\wedge\bar{\phi}=2c(a - b)\cot\alpha \, \phi\wedge\bar{\phi}, \\
&& |c|^{2} = |a|^{2} + \frac{\rho}{2}(-2 + 3\sin^{2}\alpha),
\end{eqnarray}
where (2.2) is the Gauss equation, (2.3) and (2.4) are the Codazzi-Mainardi equations, and (2.5) is the Ricci equation of $x$.
\begin{remark}  (1)\ The immersion $x$ is holomorphic (resp. anti-holomorphic) at $p \in M$ if and only if $\alpha=0\  (\mbox{resp.}\  \alpha = \pi)$ at $p$. Hence,  $\sin \alpha \neq 0$  on  $M_{0}$.

(2)\  The unitary coframes $\{w_{1},w_{2}\}$ appeared in $(2.1) - (2.5)$ are uniquely determined up to  orientations of both $\overline{M}[4\rho]$ and $M_{0}$ (see \cite{ogata}), hence the complex one-form $\phi$ on $M_{0}$ is unique up to  sign and  conjugacy. Thus, for the immersion $x$ satisfying  $(2.1) - (2.5)$, the data $\{ \psi:= \bar{\phi}, \alpha, \bar{a}, \bar{c} \}$, which satisfy  the conjugated  structure equations,  define the same surface as $x$.
\end{remark}

 For a  parallel mean curvature surface $x$ satisfying $a =\bar{a}$ on $M_{0}$, it is proved in \cite{kenzhou}, \cite{kenoga} that if $\rho >0$, then $x$ is  totally real and the image is a part of the flat torus and  in the case of  $\rho <0$, if $x$ has a constant Kaehler angle, then it  is totally real, flat and $-2b^{2} \leq \rho <0$, or it is non-totally real, $K=\mbox{constant}=-2b^{2}$, and $\rho =-3b^{2}$  (see also Theorem 1.1 of \cite{hirakawa}). If  the Kaehler angle  is not constant, then  Lemma 3.3 of \cite{kenzhou} tells us that there are `non-trivial ones'. But,  the statement of Lemma 3.3 is incomplete, so we give here the correction as follows: 
\begin{theorem} Let  $x: M  \longrightarrow
\overline{M}[4\rho]$ be  a  parallel   mean curvature surface with $|H|=2b\; (>0)$. Suppose that the Kaehler angle  is not constant and   $a=\bar{a}$ on $M_{0}$.   If $\rho <0$, then $\rho=-3 b^{2}$ and  $x$  is determined by a real constant and  the amplitude of  a Jacobi elliptic function with  the modulus  either $2\sqrt{2}/3$ or $1/3$.
\end{theorem}
Proof. All formulas and results in \cite{kenzhou} are valid under the assumption of Theorem 2.2  \cite{kenoga}.  The first claim of Theorem 2.2 is proved by the first paragraph of page 302 in   \cite{kenzhou}.  Formula (2.5)  and Theorem 4.1, both in \cite{kenzhou}, imply
%\begin{equation}
$\lambda(u)^{4}(a(u)^{2} + b^{2}(3- 9/2 \sin^{2}\alpha(u))) = m^{2}$, \
%\end{equation}
for some non-zero real constant $m$,
where $\theta(u)$ in \cite{kenzhou} is replaced here by $\alpha(u)$ here and 
$\phi =\lambda(u)(du + i dv)$.
Since $\lambda(u)$ is real-valued and defined up to a real multiplicative constant, we  assume that  $m=9 b^{3}$. Formula (2.10)  and Theorem 4.1, both in  \cite{kenzhou}, yield  $a(u) = b(1- 9/4 \sin^{2}\alpha(u))$. These two formulas imply
$\lambda(u)^{2}(8-9 \sin^{2}\alpha(u)) = \pm 36b^{2}. $
By  (2.3) of \cite{kenzhou}, $\alpha$ satisfies the ODE
%\begin{equation}
$
(d \alpha/d u)^{2} = \pm 9b^{4} (8-9 \sin^{2}\alpha(u)) .
%\end{equation}
$
Since we may assume $\alpha(u) >0$,  two cases to be studied arise:
\begin{eqnarray*}
\begin{array}{rl}
  &\mbox{ (i) } \quad  \frac{d \alpha}{d u} = 3b^{2}\sqrt{-8 +9 \sin^{2}\alpha(u)}, \quad  -8 +9 \sin^{2}\alpha(u)  >0,  \\
  &\mbox{(ii) } \quad  \frac{d \tilde{\alpha}}{d u} = 3b^{2}\sqrt{8-9 \sin^{2} \tilde{\alpha}(u)}, \qquad  8-9 \sin^{2} \tilde{\alpha}(u)  >0.
\end{array}
\end{eqnarray*}

For case (i), put $\sin \gamma(u) =-3 \cos \alpha(u)  $. The function $\gamma(u)$ satisfies the ODE
$
d \gamma/du = 9b^{2}\sqrt{1- 1/9 \sin^{2}\gamma }, \ \mbox{i.e.},
$
$\gamma(u) = am(9b^{2}u, 1/3)$, where $am(\cdot, p)$ denotes the amplitude of a  Jacobi elliptic function with the modulus $p$. Hence, 
\begin{equation}
\cos \alpha(u) = -\frac{1}{3} \sin\left (am \left(9b^{2}u,\frac{1}{3} \right)\right), \quad -\frac{K(1/3)}{9b^{2}} < u  < \frac{K(1/3)}{9b^{2}},
\end{equation}
where $K(1/3)$ denotes the complete integral of the first kind of the  Jacobi elliptic function with the modulus $1/3$.
For  case (ii), put $ \sin \tilde{\gamma}(u) = 3/(2\sqrt{2}) \sin \tilde{\alpha}(u)$. Then the function $ \tilde{\gamma}(u)$ satisfies the ODE
%\begin{equation}
$
d  \tilde{\gamma}/d u =  9b^{2}\sqrt{1 - 8/9 \sin^{2}  \tilde{\gamma} },
%\end{equation}
$
which means  $ \tilde{\gamma}(u) = am (9b^{2}u, 2\sqrt{2}/3 )$.  Hence, 
\begin{equation}
\sin \tilde{\alpha}(u) = \frac{2\sqrt{2}}{3} \sin \left (am\left(9b^{2}u, \frac{2\sqrt{2}}{3}\right) \right), \quad -\frac{K(2\sqrt{2}/3)}{9b^{2}} < u  < \frac{K(2\sqrt{2}/3)}{9b^{2}}.
\end{equation}
Now we will find the expression of $c$ in terms of $\sin \alpha$: By (2.5), $c$ is written as $c=b |8-9\sin^{2}\alpha|/4\exp(it(u,v))$ for some real-valued function $t(u,v)$. Equation (2.4) implies $t(u,v) = \mbox{constant}$.
The first and second fundamental forms of $x$ are determined by the 
Kaehler angle  $\alpha(u)$ or $ \tilde{\alpha}(u)$ and a real number $t \in [0,\pi]$ as follows:
\begin{eqnarray}
\left\{
\begin{array}{rl}
& \lambda(u) = \frac{6b}{\sqrt{|8-9 \sin^{2}\alpha(u)|}}, \quad  8-9 \sin^{2}\alpha(u) \neq 0,  \\
 & a(u) = b\left(1-\frac{9}{4} \sin^{2}\alpha(u) \right) , \quad
 c(u) = \frac{b}{4}|8-9 \sin^{2}\alpha(u)|e^{it}, 
\end{array}
\right.
\end{eqnarray}
where in the above formulas, $\alpha(u)$ is replaced by $ \tilde{\alpha}(u)$ if $8-9 \sin^{2}\alpha(u) > 0$.
Conversely,  given $b>0$,   $t \in [0,\pi]$, and the amplitude $\gamma(u)$ with a modulus $1/3$ or its conjugate  $\tilde{\gamma}(u)$, we define $\alpha(u)$ by (2.6) or (2.7), and $\lambda(u), a(u)$ and $c(u)$ by (2.8). Then it is straightforward to prove  that these data satisfy  the structure equations $(2.1)-(2.5)$ under the condition $\rho =-3b^{2}$,  proving Theorem 2.2.
\begin{remark}
Case (i)  above is missed in Lemma 3.3 of  \cite{kenzhou}.
  \end{remark}
 We denote both surfaces obtained  in Theorem 2.2 as $x_{t,1/4+}$ if it is determined by $\gamma(u)$ and $x_{t,1/4-}$ if  determined by $\tilde{\gamma}(u)$.   The surface $x_{t,1/4\pm}$ is isometric to  $x_{0,1/4\pm}$ and  both surfaces  have the same length of the mean curvature vector.  Consequently,  the set  $\{x_{t,1/4\pm}\; |\; t \in [0,\pi] \}$ is the associated family of $x_{0,1/4\pm}$. We note that in general $\sin \alpha(u) \neq \sin \tilde{\alpha}(u)$, hence  surfaces $x_{t,1/4+}$ and $x_{t,1/4-}$ are not isometric to each other.
Let $\mathbb{D}_{\pm}$ be the domain in $\mathbb{R}^{2}$ given by
$$
\mathbb{D}_{\pm} = \left\{ (u,v) \in \mathbb{R}^{2} \ | -\frac{K(q)}{9b^{2}} < u  < \frac{K(q)}{9b^{2}}, \  -\infty < v < \infty \right\},
$$
where the plus  sign and the minus one of $\mathbb{D}_{\pm}$ are taken for $q=1/3$ and  $q=2\sqrt{2}/3$, respectively.
 We will prove that the immersions $x_{t,1/4\pm}$ are complete.   
  \begin{proposition}  For each $t \in [0,\pi]$, 
     $x_{t,1/4\pm} :  \mathbb{D}_{\pm} \longrightarrow   \overline{M}[-12b^2]$ are complete.
   \end{proposition}
   Proof. For case (i),  the Riemannian metric $ds^{2}$ on $x_{t,1/4+}$ is estimated as follows:
   \begin{eqnarray*}
 \frac{1}{36b^{2}} ds^2 &=&  \frac{du^2 + dv^2}{-8+9 \sin^{2}\alpha} = \frac{du^2 +dv^2}{\cos^{2}\gamma(u)}, \quad  \left(-\frac{\pi}{2} < \gamma < \frac{\pi}{2} \right) \\
  &=& \frac{1}{\cos^2 \gamma} \left( \frac{d\gamma^2}{9 (1- 1/9 \sin^2 \gamma)} + dv^2 \right) \geq \frac{1}{9} \frac{1}{\cos^2 \gamma} \left(d\gamma^2 + d(3v)^2 \right).
   \end{eqnarray*}
Since the metric in the right hand  is complete ( see, for instance, Lemma 3.13 of  \cite{hirakawa}), the  considered metric is  also complete.  The proof in the second case is  similar and we conclude the proof.

\begin{remark}
In  the third line  of the proof of  Lemma 3.13 in  \cite{hirakawa},  $\arctan u$ should be read as $\tan u$.
\end{remark}

A parallel mean curvature surface  is of {\it a general type} if it satisfies $ d\alpha \neq 0$ and  $a \neq \bar{a}$  at a point of $M_{0}$.
Now, we will  recall  results in \cite{ken1} and \cite{ken2}   which are used in the next section:  if $x$ is of a general type, then
$a$ is a function of $\alpha$, say $a=a(\alpha)$,  satisfying the first order complex ODE (see (3.1) of \cite{ken1})
\begin{equation}
 \frac{da}{d\alpha} = \frac{\cot \alpha}{\overline{a} +b}\left( - 2ba +  2|a|^{2}
 + \frac{3\rho}{2}\sin^{2}\alpha \right), \ (a+b \neq 0).
 \end{equation}

Let us consider  isothermal coordinates for the Riemannian metric $\phi \bar{\phi}$ on $M_{0}$ that makes $M_{0}$ a Riemann surface with a local complex coordinate $z$.
%The following Lemma was proved in \cite{ken2}
%Suppose that $x$ is of  a general type.  
By Lemma 2.4 of \cite{ken2}, there is a complex analytic transformation $w=w(z)$ on the Riemann surface $M_{0}$ such that  in $\phi = \mu dw$, $\mu$ is a complex-valued function of  a single real variable  $(w+\bar{w})/2$.  Let $w=u+i v, \  (u, v \in \mathbb{R})$.  Then $\alpha$ is a function of $u$ only, satisfying the ODE
\begin{equation}
\frac{d\alpha}{du} = 2 \exp \left(\int F(\alpha)d\alpha \right),
\end{equation}
where
%\begin{equation}
$F(\alpha) =\cot \alpha (|a - b|^{2} + 3\rho/2 \sin^{2} \alpha)/|a + b|^{2}$  
%\end{equation}
%$\alpha_{z} = G(z) \exp(\int F(\alpha)d\alpha)$. We set $w= \int G(z)dz$. 
%$\phi = \lambda dz = \alpha_{z}/(a(\alpha)+b)dz=\mu dw,$
and $\mu$ is expressed as
\begin{equation}
\mu =  \frac{\exp \left(\int F(\alpha)d\alpha \right)}{a+b},
\end{equation}
where  $\mu$ is a  complex-valued function of $\alpha$  defined up to sign and a real multiplicative  constant.

We will also need  the following fact  proved in  Lemma 2.5 of \cite{ken2}:
there are  real numbers $k$ and $t$ such   that
\begin{equation}
c= \left(|a|^{2} + \frac{\rho}{2}\left(-2 + 3 \sin^{2}\alpha \right) \right)^{1/2} \frac{(\overline{a}+b)}{a+b} e^{-i (k v+ t)}.
\end{equation}
\begin{remark}
In (2.14) of \cite{ken2}, $exp(-it), \; (0 \leq  t \leq \pi)$, is missed,  because $c$ is unique up to  a complex multiplicative constant with the absolute value equal to one.
\end{remark}

\section{Main lemma}
In this section, we prove that the constant $k$ in $(2.12)$ is zero, which is the key result of this article.

\begin{lemma} Let  $x: M  \longrightarrow
\overline{M}[4\rho]$ be a parallel mean curvature surface with $|H|=2b\; (>0)$ of a general type. Then
\begin{eqnarray}
|a +b|^{2} &=&  \frac{(8b^{2} + 3\rho \sin^{2}\alpha)}{4 b} ( \Re a +  b) ,  \\
(\Im a)^{2} &=& (\Re a +b) \left(\frac{(8b^{2}+3\rho\sin^{2}\alpha)}{4b} - (\Re a +b) \right).
\end{eqnarray}
where $\Re a$ and $\Im a$ denote the real  and  imaginary parts of $a$ respectively.
\end{lemma}
Proof.  Formula  (2.17) of  \cite{ken2} yields
$
\mu^{2}(8ba - 3 \rho \sin^{2}\alpha) = c_{1},
$
where the constant $c_{1} (\in \mathbb{C})$  is not zero, because the immersion $x$ is of a general type. By a transformation: $\mu \rightarrow   e^{i\xi} \mu,\; (\xi \in \mathbb{R})$,   $c_{1}$ can be assumed to be real. By (2.11), $\mu \cdot (a+b)$ is also real-valued. Coupling these together,  one obtains
$$
\mu^{2}\cdot(a+b)^{2} = \frac{c_{1}(a+b)^2}{8ba - 3 \rho \sin^{2}\alpha} =\overline{\mu^{2}\cdot (a+b)^{2}} = \frac{c_{1}(\bar{a}+b)^{2}}{8 b \bar{a} -3 \rho \sin^{2}\alpha}.
$$
%which yields
%$(a+b)^{2}(8 b \bar{a} -3 \rho \sin^{2}\alpha)=(\bar{a}+b)(8 b a -3 \rho \sin^{2}\alpha).
%$
The simplification of the above formula  proves  (3.1). Formula (3.2) is another expression of (3.1), proving Lemma 3.1.

Since we are now studying a parallel mean curvature surface of a general type,  $\alpha$ is not constant, hence $8 b^{2} + 3 \rho \sin^{2}\alpha \neq 0$. Then the real part of $a$ satisfies a linear ODE:
\begin{lemma} Let  $x: M  \longrightarrow
\overline{M}[4\rho]$  be a parallel mean curvature surface of a general type. Then $\Re a$  satisfies  the first order linear ODE
\begin{equation}
 \frac{d}{d\alpha}\Re a+ 2 \cot \alpha \frac{ (8 b^{2}  - 3\rho \sin^{2}\alpha )}{8 b^{2}+3 \rho \sin^{2}\alpha }  \Re a - 2b\cot \alpha  \frac{(8 b^{2} + 9\rho \sin^{2}\alpha) }{8 b^{2}+3 \rho \sin^{2}\alpha} =0 .
\end{equation}
\end{lemma}
Proof. 
It follows from (2.9) and (3.2) that
$$
 \frac{d}{d\alpha}\Re a + \frac{(\Re a+b)}{2b|a+b|^{2}} \left\{(8b^{2} - 3 \rho \sin^{2}\alpha)\Re a - b(8b^2 + 9\rho \sin^{2}\alpha ) \right\} \cot \alpha =0.
$$
We conclude by coupling this with (3.1).

 If $\rho \neq 0$, then the general solution of (3.3) is given by
\begin{equation}
\Re a = -b + \frac{1}{3b\rho \sin^{2}\alpha}  (8b^{2}+3\rho \sin^{2}\alpha)( - 2b^{2} +c_{3}(8b^{2} + 3 \rho \sin^{2}\alpha)), \quad (c_{3} \in \mathbb{R}),
\end{equation}
and if $\rho =0$,  then
$\Re a = (c_{3}- b/2 \cos 2 \alpha)/\sin^{2}\alpha$.

By (3.1) and (3.4),  
\begin{equation}
|a+b|^{2} = \frac{(8b^{2}+3\rho \sin^{2}\alpha)^{2}}{12b^{2}\rho\sin^{2}\alpha} (-2b^{2} + c_{3}(8b^{2}+3\rho \sin^{2}\alpha)), \quad (\rho \neq 0),
\end{equation}
in particular,
\begin{equation}
\rho (-2b^{2} + c_{3}(8b^{2}+3\rho \sin^{2}\alpha) \geq 0.
\end{equation}
The formula $|a-b|^{2}  = |a+b|^{2} - 4 b \Re a$, (3.1), and (3.4) yield
%\begin{equation}
$$
 |a-b|^{2} + \frac{3\rho}{2}\sin^{2}\alpha 
 = \frac{(8b^{2}+3\rho\sin^{2}\alpha)(16 b^{4} - c_{3}(64 b^{4}-9 \rho^{2} \sin^{4}\alpha)) }{12 b^{2} \rho \sin^{2}\alpha}. 
%\end{equation}
$$
By (3.4) and the above formula, the function $F(\alpha)$  in (2.10) can be expressed as
$$
F(\alpha) =  \frac{(16 b^{4} - 64 b^{2} c_{3} + 9 c_{3} \rho^{2} \sin^{4}\alpha)}{(8 b^{2} + 3 \rho \sin^{2}\alpha)(-2 b^{2} + c_{3} ( 8 b^{2} + 3 \rho \sin^{2}\alpha))}\cot \alpha.
$$
Integrating, we obtain, up to a positive multiplicative constant,
\begin{equation}
e^{\int F(\alpha)d\alpha} =\frac{\sqrt{|8b^{2}+3\rho \sin^{2}\alpha||-2 b^{2} + c_{3}(8b^{2}+3\rho \sin^{2}\alpha)|}}{\sin \alpha}.
\end{equation}
By (3.2) and (3.4), the imaginary part of $a$  is expressed 
 in terms of $\sin \alpha$ as follows:
\begin{equation}  
(\Im a)^{2} = \frac{(1-4 c_{3})}{36 b^{2} \rho^{2}\sin^{4}\alpha} (8b^{2}+3\rho \sin^{2}\alpha)^{3} (-2b^{2} + c_{3}(8 b^{2}+3\rho \sin^{2}\alpha)), \quad  (\rho \neq 0).
\end{equation}

\begin{remark}
We may assume $\Im a \geq 0$, by Remark 2.1 (2).
\end{remark}
If $\rho \neq 0$, Equation (3.8)  implies
\begin{equation}
 (1-4c_{3})(8b^{2}+3\rho \sin^{2}\alpha)(-2b^{2} + c_{3}(8 b^{2}+3\rho \sin^{2}\alpha)) \geq 0.
\end{equation}
%\begin{remark} When $\rho =-3b^{2}$, (3.8) is identical to (3.9) of \cite{ken2}.
%\end{remark} 
The next lemma is the main result of this section.
\begin{lemma} Let  $x: M  \longrightarrow \overline{M}[4\rho]$ be a parallel mean curvature surface of a general type. If $\rho \neq 0$, then  $\rho = - 3 b^{2} < 0$ and $k=0$ or $c_{3}=c \equiv 0$. Moreover, in the latter case  the Gaussian curvature of $M$ is 
constant  equal to $=-2b^{2}$.
\end{lemma}
Proof. 
By (2.16) of \cite{ken2} and (2.11) of this article, we have:
%\begin{equation}
%$$2k
%= \rho\mu \cdot  (a+b)\frac{\left( 8 |a|^2 + 9b(a+\bar{a}) \sin^2 \alpha - 8 b^2 + 18 b^2\sin^2\alpha\right)}{ |a+b|^{2} |c|^{2}}  \cot \alpha
%\end{equation}
%$$
%which implies, by (2.12), 
\begin{equation}
2k|a+b|^{2}|c|^{2}
= \rho e^{\int F(\alpha)d\alpha}( 8 |a|^2 +  9 b (a+\bar{a})\sin^2 \alpha - 8 b^2 + 18 b^2)\cot\alpha ,
\end{equation}
where $k_{1}$ in \cite{ken2} is replaced here by  $k$.
It follows from (2.5),  (3.1) and (3.4) that
\begin{equation}
%\begin{displaymath}
\left\{ \begin{array}{l}
  8 |a|^2 + 9b(a+\bar{a}) \sin^2 \alpha - 8 b^2 + 18 b^2\sin^2\alpha  
 =\frac{6(\rho + 3b^{2})}{b}(\Re a+b) \sin^{2}\alpha, \\
|c|^{2}=  -(\rho + 3b^{2}) + \frac{c_{3}}{4b^{2}} (8b^{2} + 3\rho \sin^{2}\alpha)^{2}.  
\end{array}  \right.
%\end{displaymath}
\end{equation}

By (3.5),  (3.7), (3.10), and the above formulas, we obtain
\begin{eqnarray}
&& \rho^{2}(\rho+3 b^{2})^{2}(-2 b^{2}+c_{3}(8b^{2}+3\rho \sin^{2}\alpha)  ) 
 (1-\sin^{2}\alpha) \\
&&  \qquad  \qquad  \qquad  - k^{2}(8b^{2}+3\rho \sin^{2}\alpha)\left(  -(\rho + 3b^{2}) + \frac{c_{3}}{4b^{2}} (8b^{2} + 3\rho \sin^{2}\alpha)^{2}\right)^{2} =0.  \nonumber 
\end{eqnarray}
This is a polynomial equation of degree 5 in $\sin^{2}\alpha$ with constant coefficients. It follows that if $\rho \neq 0$ and $\alpha$ is not constant, then all the coefficients of the polynomial vanish, in particular,  $k=0$ or $c_{3}=0$ which comes from   the term of the  highest degree.  Then we have  $\rho =- 3 b^{2}$ by (3.12).  In the case when $c_{3}=0$,  we get $c \equiv 0$ by (3.11) and  $K=\mbox{constant} = -2b^{2}$ by  (2.2) and (2.5).  Consequently, this is Hirakawa surface \cite{hirakawa}, which concludes the proof.

It shall be noted that Lemma 3.4 does not need any global assumption of $M$.
%By Lemma 4.4,  all results in section 2 of \cite{ken2} hold without the assumption $k=0$. 
As a direct application of Lemma 3.4,  we have:
%\begin{theorem}
{\it  Even locally, there is  no parallel mean curvature surface of a general type   in $\mathbb{C}P^{2}$.}
%\end{theorem}
 So, a parallel mean curvature surface in $\mathbb{C}P^{2}$ has either constant Kaehler angle or non-constant Kaehler angle satisfying $a=\bar{a}$ on an open subset of $M$. However, the latter does not occur by Theorem 4.1 of \cite{kenzhou}.  (We note that all results in \cite{kenzhou} hold under the additional condition $a=\bar{a}$  \cite{kenoga}. ) Therefore, a parallel mean curvature surface in $\mathbb{C}P^{2}$ must have a constant Kaehler angle, in particular, a constant Gaussian curvature.  Then, such a surface is flat and totally real \cite{hirakawa}. Thus, we proved
\begin{theorem}
Any complete parallel mean curvature surface with $|H|=2b\; (>0)$ in $\mathbb{C}P^{2}$ is  the flat torus up to  rigid motions of $\mathbb{C}P^{2}$.
\end{theorem}
Theorem 3.5 is in contrast to the fruitful theory of minimal surfaces in $\mathbb{C}P^{2}$ (see Introduction of  \cite{inotaniudag} for references). 

\section{Jacobi elliptic functions}
In this section, non-zero parallel mean curvature surfaces   in $\mathbb{C}H^{2}$  are studied by applying the results of Section 3,  and we will prove that such a  surface  is determined by  a Jacobi elliptic function.   
Let $x: M \longrightarrow \bar{M}[4\rho]$ be a parallel mean curvature surface of a general type with $|H|=2b\; (>0)$. We  can assume $\rho = - 3b^{2}$ by Lemma 3.4 and $a \neq \bar{a}$ by Theorem 2.2. Lemma 3.4 and the second formula of  (3.11) yield  $c_{3} \geq 0$, and
$a$ is explicitly expressed in terms of  $\sin \alpha$ by (3.4),  (3.8), and Remark 3.3:
\begin{eqnarray}
a &=&\frac{b}{9\sin^{2}\alpha}\{ 16 - 27 \sin^{2}\alpha - c_{3}(8 - 9 \sin^{2}\alpha)^{2} \nonumber \\
&&  \qquad  +  \frac{i}{2}|8  - 9 \sin^{2}\alpha|\sqrt{(1-4c_{3})(8-9 \sin^{2}\alpha)(-2 + c_{3}(8  - 9 \sin^{2}\alpha))}\  \}.
\end{eqnarray}
From now on, let $p=c_{3}$. By  (3.6) with $\rho + 3 b^{2}=0$ and (3.9),  since $2-p (8- 9 \sin^{2}\alpha) >0$, two cases to be studied arise:
%\begin{itemize}
%\item 
  (i) \ $0<p < 1/4$ and $ 8 - 9 \sin^{2}\alpha <0$, and
 (ii) \ $1/4 < p$ and $8-9\sin^{2}\alpha >0$.
%\end{itemize}

First, we study case (i).  By (2.10) and (3.7),  $\alpha=\alpha(u)$ satisfies the ODE
\begin{equation}
\frac{d\alpha}{du} = 2b^{2}\frac{\sqrt{(-8 +9\sin^{2}\alpha)(2-p (8- 9 \sin^{2}\alpha))}}{\sin \alpha}.
\end{equation}
%if necessary, $u$ is changed to $-u$.
Consider the function  $\gamma = \sin^{-1}(-3 \cos \alpha)$. By (4.2), this satisfies the ODE
 \ $d\gamma/du = 6b^{2}\sqrt{2+p - p \sin^{2}\gamma}$.
Hence, $\gamma$ is the amplitude of  a Jacobi elliptic  function with the modulus $\sqrt{p/(2+p)}$, i.e.,  $\gamma = \gamma(u)$  is written as   
\begin{equation}
\gamma = am(6b^{2} \sqrt{2+p}u, \sqrt{p/(2+p)}) , \quad  -\frac{\pi}{2} < \gamma < \frac{\pi}{2}.
\end{equation}
For case (ii),  $\alpha=\alpha(u)$ satisfies the ODE
\begin{equation}
\frac{d\alpha}{du} = 2b^{2}\frac{\sqrt{(8 -9\sin^{2}\alpha)(2-p (8- 9 \sin^{2}\alpha))}}{\sin \alpha}.
\end{equation}
Let $\tilde{\gamma} = \cos^{-1} \sqrt{p(8-9 \sin^{2}\alpha)/2}$.
This satisfies the ODE \  $d\tilde{\gamma}/du = 6b^{2}\sqrt{2+p - 2 \sin^{2}\tilde{\gamma}}$.
Hence, $\tilde{\gamma}$ is the amplitude of  a Jacobi elliptic  function with the modulus $\sqrt{2/(2+p)}$, i.e.,   $\tilde{\gamma} = \tilde{\gamma}(u)$  is written as   
\begin{equation}
\tilde{\gamma}  = am(6b^{2} \sqrt{2+p}u, \sqrt{2/(2+p)}), \quad  -\frac{\pi}{2} <\tilde{\gamma} < \frac{\pi}{2}.
\end{equation}
We note that $\tilde{\gamma} $ is the conjugate of $\gamma$.
We showed that in  these two cases, {\it  Kaehler angles are  determined by a Jacobi elliptic function or its conjugate}.

Let us put
 \begin{equation}
\mathbb{D}= \left\{(u,v) \in R^{2}\  |\  -\frac{K(q)}{6b^{2}\sqrt{2+p}}<u <\frac{K(q)}{6b^{2}\sqrt{2+p}},\   -\infty <v< \infty \right\},
\end{equation}
where  $K(q)$ denotes the complete elliptic integral of the first kind with the modulus $q=\sqrt{p/(2+p)}$ for $0<p<1/4$ and $q=\sqrt{2/(2+p)}$ for $p>1/4$. We note that $\mathbb{D}$ is transformed to the slab $\{ (\tilde{u}, v) \in \mathbb{R}^{2}\  | -\pi/2 < \tilde{u} < \pi/2, \ v \in \mathbb{R} \}$, where $\tilde{u} =\gamma$ for $0<p<1/4$ and $\tilde{u}= \tilde{\gamma}$ for $p>1/4$. 
%$$
%8 - 9 \sin^{2}\alpha >0, \quad 2 - p(8-9 \sin^{2}\alpha) >0, \quad (\frac{1}{4} < p).
%$$

Given $p \in (0,\infty)$, 
let us define a real-valued function of $\alpha$, $\theta(\alpha)$, by
\begin{equation}
e^{i\theta(\alpha)} =  \frac{1}{3\sin \alpha}\left( 2\sqrt{2-p(8-9\sin^{2}\alpha)} 
\pm i \sqrt{|(1-4p)(8-9\sin^{2}\alpha)|} \right),
\end{equation}
where the plus sign is taken for $0 < p  \leq 1/4$ and the minus one for $1/4 < p$.
By (2.11), (2.12), (3.7) and (4.1),  $\mu, a,$ and $c$  are expressed as follows:
\begin{eqnarray}
%\left\{
%\begin{array}{rl}
\mu &=&   \frac{6b e^{i\theta(\alpha)}}{ \sqrt{|8 -9 \sin^{2}\alpha|} },\quad  (8 -9 \sin^{2}\alpha \neq 0 ),  \nonumber \\
a &=& -b + \frac{b(8-9\sin^{2}\alpha)\sqrt{2- p(8 -9 \sin^{2}\alpha)}}{6\sin \alpha} e^{-i\theta(\alpha)} \\
c &=& \frac{b}{2}\sqrt{p}|8-9 \sin^{2}\alpha| e^{i(2\theta(\alpha) +t)}, \ (t \in  [0,\pi]). \nonumber
%\end{array}
%\right.
\end{eqnarray}
The last main result of this article is the following theorem.
\begin{theorem} Let  $x: M  \longrightarrow
\overline{M}[-12b^{2}]$ be  a parallel   mean curvature surface  with $|H|=2b\;(>0)$. If $x$ is of a general type, then there exists a positive number  $p\; (\neq 1/4)$ such that the Kaehler angle  is determined by a Jacobi elliptic function $(4.3)$ if $0 < p < 1/4$ and $(4.5)$  if $1/4 < p  <\infty$. Moreover,  the first and  second fundamental forms of $x$ are given by $(4.8)$. Conversely, given  three real numbers $b\; (>0) , t\in [0,\pi],\  p\; (>0, \neq 1/4)$, there exist a domain $\mathbb{D} \subset \mathbb{R}^{2}$ and a parallel mean curvature immersion of a general type  from $\mathbb{D}$ into $\overline{M}[-12b^{2}]$ such that the first and  second fundamental forms are determined by  $t$ and a Jacobi elliptic function whose  modulus depends only on $p$, and the length of the mean curvature vector is equal to $2b$.
\end{theorem}
Proof. We only  have to prove the last part of Theorem 4.1.   Consider three numbers $b,t, p$ presented in Theorem 4.1 and define the domain $\mathbb{D} \subset \mathbb{R}^{2}$ by (4.6). Suppose, first, that $0 < p < 1/4$. Take  the amplitude  $\gamma$ of  a Jacobi elliptic  function as in (4.3).
%$\gamma = am(6b^{2} \sqrt{2+p}u, \sqrt{p/(2+p)}) .$ 
Define  a function $\alpha$ by 
$\alpha =\alpha(u; p) = \cos^{-1}(-\sin \gamma/3 )$.
Then, 
$8-9\sin^{2}\alpha <0, \ 2-p(8-9\sin^{2}\alpha)>0$,
and  $\alpha$ satisfies  the ODE  (4.2).
If $p> 1/4$,  for a given amplitude $\tilde{\gamma}$  as in  (4.5), we have a function $\alpha$ which satisfies  the ODE  (4.4).
Let us  define a real-valued function $\theta(\alpha)$ by (4.7), and $\mu, a, c$ by (4.8). Now, we prove that
%\begin{proposition}
the data $\{\mathbb{D}, \phi = \mu (du+idv), a,c \}$  satisfy  the structure equations $(2.1) - (2.5)$ under the condition  $\rho = - 3b^{2}$. 
%\end{proposition}
In fact,   by the first and  second formulas of (4.8),  we see that $ (a+b) \mu = (\bar{a}+b) \bar{\mu}$,
which proves the first equation of (2.1). To show the second equation of (2.1), we 
use that an equation $ d\phi 
=-\log \mu'(\alpha) (\bar{a} +b)   \phi \wedge\bar{\phi}$  which is derived by virtue of the first equation of (2.1). Thus, the second equation of  (2.1)  is equivalent to
$
\mu'(\alpha) (\bar{a}+b) + \mu (\bar{a}-b) \cot \alpha =0.
$
This follows from (4.8) by direct computation.  
To prove (2.2), we recall the formula of  Gauss curvature $K =- e^{-2\sigma}(\sigma_{uu} + \sigma_{vv})$ in terms of isothermal coordinates $(u,v)$ given by $ds^{2} =e^{2\sigma}(du^2  + dv^2)$. Now,  $e^{2\sigma} = 36 b^2/(-8+9 \sin^{2}\alpha(u;p))$, and  the direct computation using (4.2)  implies
$
\sigma_{uu}= 36 b^{4}(2+ p (-8+9 \sin^{2}\alpha)^2)/(-8+9 \sin^{2}\alpha),
$
 which implies  (2.2). To prove (2.3), we  first note that
$
\mu^{2}(8 a + 9b \sin^{2}\alpha) =  \mbox{constant},
$  which is observed by (4.8).
Taking the exterior derivative of this formula, (2.3) with $\rho=-3 b^2$  follows from (2.1).
To prove (2.4), we use
$\mu^2 \bar{c} = \mbox{constant}$  which can be seen from (4.8). Taking the exterior derivative of this formula and using (2.1), (2.4) is proved.
 (2.5) follows from  the formula $|a|^{2}=  |a+b|^2 - 2b \Re (a+b)  + b^2$  and  (4.8). For the case when $p > 1/4$, the proof is done by the similar way as before  using (4.4) and (4.5), so we omit it, proving  Theorem 4.1.

Let $\mathbb{C}H^{2}[4\rho]$ be the complex hyperbolic plane  of constant holomorphic sectional curvature $4\rho,\ (\rho <0)$.
By $x_{t,p}, \ (0 \leq t \leq  \pi, \ 0<p<\infty, \ p \neq 1/4)$, we denote the  immersion from $\mathbb{D} \subset \mathbb{R}^{2}$ into  $\mathbb{C}H^{2}[-12b^{2}]$ obtained by 
Theorem 4.1. 
 Note that these surfaces are the simplified version of  the ones given in Section 4 of  \cite{ken1}, because $\xi(t)$ in \cite{ken1} is explicitly integrated here by using the special isothermal coordinates. %So, this expression in this paper  is much simpler than those of the formulas in \cite{ken1}.

Now, we list the  properties of  $x_{t,p}$:
\begin{enumerate}
\item The  mean curvature vector $H$ is parallel,  $|H|=2b$, and  the Kaehler angle  is not constant.

It follows from the  structure equations $(2.1)- (2.5)$.
\item  $x_{t,p} : \mathbb{D} \longrightarrow \mathbb{C}H^{2}[-12b^{2}]$ is complete. 
 
We do not give the proof here, as it can be performed almost  in the same way as the proof of Proposition 2.4.

\item For any $p,q$ with $0<p<q<\infty$, $x_{t,p}$ is not isometric to $x_{t,q}$, and hence they are not congruent each other, because of $\sin^{2}\alpha(u; p) \neq \sin^{2}\alpha(u; q)$.
\item   By letting $p \rightarrow 0$, $x_{t,p}$ has a limit, say $x_{t,0}$,  that is  Hirakawa surface \cite{hirakawa}. 

In fact, as $p \rightarrow 0$ in (4.2), (4.7) and (4.8),  $\alpha$ satisfies an ODE, 
%$$
%\frac{d\alpha}{du} = 2\sqrt{2}b^{2}\frac{\sqrt{-8 + 9 \sin^{2}\alpha}}{\sin \alpha}, 
%$$ 
$\phi$ and $a$ are well defined and $c \equiv 0$.  It is verified in the same way as in the case when $0 < p < 1/4$ in  Theorem 4.1 that these data satisfy the structure equations $(2.1) - (2.5)$ under the condition $\rho = -3 b^{2}$. Therefore, there exists an immersion $x_{t,0} : \mathbb{D} \longrightarrow \mathbb{C}H^{2}[-12b^{2}]$. It can be easily shown that the immersion $x_{t,0}$ has a constant Gaussian curvature $K=-2b^{2}$ by (2.2), but is not of a constant Kaehler angle. So, $x_{t,0}$ is congruent to Hirakawa surface by Theorem 1.1 (1) of  \cite{hirakawa}.

\begin{remark}
 At page 230 of   \cite{hirakawa}, author describes two parallel mean curvature surfaces of constant Gaussian curvature $K=-2b^{2}$, but both surfaces are congruent to each other in $\mathbb{C}H^{2}[-12b^{2}]$ by Remark 2.1 (2).
\end{remark}

 \item As  $p \rightarrow 1/4-0$, $x_{t,p}$ has  a limit which is $x_{t,1/4-}$,  and as  $ p \rightarrow 1/4+0$, $x_{t,p}$ has also a limit which is $x_{t,1/4+}$.  

In fact, as  $p \rightarrow 1/4-0$ in (4.2),  (4.7), and (4.8), $\alpha, \phi, a$ and $c$ have  limits which are the same ones in Theorem 2.2.   In the  case when $p>1/4$, we have  similar formulas.
%These are obtained in  Kenmotsu and Zhou \cite{kenzhou}.
\item  $x_{t,p}$ is isometric to $x_{\tilde{t},p}$ for any $t, \tilde{t} \in [0,\pi]$  and  both surfaces have the same length of the  mean curvature vector. Hence, for any fixed $p$,  $\{x_{t,p}  \ |\  t \in [0,\pi] \}$ is the  associated family of  $x_{0,p}$.

In fact,  the metric is determined only by  the Kaehler angle function and thereofore independent from $c$.
 \end{enumerate}

Finally,  we describe the manifold structure of  the  set of the complete  parallel mean curvature surfaces in $\mathbb{C}H^{2}[4\rho]$ as follows:
By {\it the moduli space of parallel mean curvature surfaces}, we mean the quotient set of all parallel mean curvature surfaces, where the equivalence relation is given by  rigid  motions of the target space. By Theorem 3.5, the moduli space of complete  parallel mean curvature surfaces in $\mathbb{C}P^{2}$ reduces to a point.

 If there exists a parallel mean curvature surface of a constant Gaussian curvature in  $\mathbb{C}H^{2}[4\rho]$, then $\rho$ satisfies the condition  (i) $-2b^{2} \leq \rho < 0$ or (ii) $\rho = -3b^{2}$  \cite{hirakawa}. For case (i), by Theorem 2.2 and Lemma 3.4, a parallel mean curvature surface with $|H|=2b$ in $\mathbb{C}H^{2}[4\rho], \ (-2b^{2}\leq \rho <0),$ is of a  constant Kaehler angle, and hence has a  constant Gaussian curvature. Then, the surface is flat \cite{hirakawa} and the image is explicitly determined in \cite{hirakawa2}.
For case (ii), if the   Kaehler angle  is constant,  then the surface is uniquely determined as  Chen surface \cite{chen3} and if the Kaehler angle  is not constant,  then  we have Hirakawa surface \cite{hirakawa}. Applying Theorem 2.2 and Theorem 4.1, we get  the following:
\begin{theorem}
The moduli space of complete parallel mean curvature surfaces with $|H|=2b$ in 
$\mathbb{C}H^{2}[-12b^{2}]$ is the disjoint union of   two cones in $\mathbb{R}^{3}$.
\end{theorem}
Proof. 
Let $x: M \longrightarrow \mathbb{C}H^{2}[-12b^{2}]$ be  a   parallel  mean curvature surface with $|H|=2b>0$, where it is not assumed that $x$ is of a general type.  First, we  prove that if {\it $a \neq \bar{a}$ at a point of $M_{0}$, then this holds throughout $M_{0}$}. In fact, suppose that  there is a point ${\bold u}_{0} \in M_{0}$ with $a =\bar{a}$ at ${\bold u}_{0}$.  Since the set of points in  $M_{0}$  satisfying $a=\bar{a}$ is not open,  as $n \rightarrow \infty$,  we have a sequence of points ${\bold u}_{n} \rightarrow {\bold u}_{0}$ such that  $a \neq \bar{a}$ at  ${\bold u}_{n}$. Then, at ${\bold u}_{n}$, we have (4.1) for $a$  and, in particular, $1-4p \neq 0$.  As $n$ tends to $\infty$,  both functions  $8-9 \sin^{2}\alpha$ and $2-p(8-9 \sin^{2}\alpha)$ do not vanish at ${\bold u}_{0}$, because, otherwise, $a \rightarrow -b \ ( {\bold u}_{n}  \rightarrow {\bold u}_{0})$   by the second formula of (4.8), which contradicts $|\mu| < \infty$ at ${\bold u}_{0}$ by (2.11). Consequently,   since the imaginary part in (4.1) goes to zero, we have  $1-4p =0$,  giving a contradiction.  We proved that {\it the set of parallel mean curvature surfaces in $\mathbb{C}H^{2}[-12b^{2}]$ is divided by the surfaces with $a=\bar{a}$ on $M_{0}$ and
 those with $a \neq \bar{a}$ on $M_{0}$.}
If  a complete parallel mean curvature surface with $|H|=2b$  in $\mathbb{C}H^{2}[-12b^{2}]$  is of a general type,  then there are  real numbers $p\; (>0)$  and $t \in [0,\pi]$ such that $x = x_{t,p}, \ p \neq 1/4$, up to  rigid motions of $\mathbb{C}H^{2}[-12b^{2}]$ by Theorem 4.1, i.e., $(t,p)$ are the coordinates of the non-singular part of the moduli space we now considering.   Since $x_{0,p}$ and $x_{\pi,p}$ are congruent because $c$ is defined up to  sign, we  identify $t=0$ and $t=\pi$ in $[0,\pi]$, so we have a circle $\mathbb{S}^{1}$, and the moduli space  of the surfaces satisfying $d\alpha \neq 0$ and $a \neq \bar{a}$ is the disjoint union  of $\mathbb{S}^{1} \times (0,1/4) \cup \mathbb{S}^{1} \times (1/4, \infty)$. When $p=0$, up to  rigid motions, there is only one surface of $\mathbb{C}H^{2}[-12b^{2}]$ which is  Hirakawa surface,  and as $p\rightarrow 1/4-$, there is the one-parameter  family  $\{x_{t,1/4-}\ | \  t  \in \mathbb{S}^{1}\}$. These make a $Cone_{1}$ in $R^{3}$, in which the vertex of the cone is    Hirakawa surface and the base of the cone is  the circle  $\{x_{t,1/4-} \ | \ t \in [0,\pi] \}$.   A $Cone_{2}$ is  the one point compactification of infinite cylinder 
${\mathbb S}^{1} \times [1/4, \infty)$, where the base  of the cone is the circle $\{x_{t,1/4+}
\ | \ t  \in {\mathbb S}^{1} \}$ and  $\infty$ corresponds to Chen surface,
 proving Theorem 4.3.

  Compact parallel mean curvature surfaces with the genus $\leq 1$ in $\mathbb{C}P^{2}$ and $\mathbb{C}H^{2}$ are classified in \cite{ogata},   \cite{hirakawa},  \cite{fetcu},  \cite{ken2}.  As a corollary of   Theorem 4.3 of this article, we can see that the assumption of the genus in those classification results is not necessary.
\begin{corollary} Let $M$ be a two-dimensional compact Riemannian manifold and  let $x: M \longrightarrow \overline{M}[4\rho]$  be an isometric immersion with a non-zero  parallel  mean curvature vector.  If $\rho \neq 0$, then both  the Kaehler angle  and the Gaussian curvature are constant.  
\end{corollary} 
Parallel mean curvature surfaces of constant Gaussian curvature are even locally classified in  \cite{hirakawa}.

\section{Added in proof}
We can state two open problems for parallel mean curvature surfaces which are related to this article.
\begin{enumerate}
\item    A parallel mean curvature surface in a higher dimensional complex space-form is contained in a totally geodesic submanifold of the ambient space with  codimension at most 8, if the surface is not totally real \cite{fetcu}, \cite{fertri}.   Can be the main results of this article, Theorem 3.5 and Theorem 4.1,   generalized to  $\mathbb{C}P^{n}$ and $\mathbb{C}H^{n}$ with   $ (3 \leq n \leq 5)$ ? 
\item   Parallel mean curvature surfaces in $\mathbb{S}^{3} \times \mathbb{R},\; \mathbb{H}^{3} \times \mathbb{R},\; \mathbb{S}^{2} \times \mathbb{S}^{2}$, and $\mathbb{H}^{2} \times \mathbb{H}^{2}$ are studied in \cite{fetros1} and \cite{torurb}. What are the  moduli spaces of  the complete parallel mean curvature surfaces in these ambient spaces ?
\end{enumerate}

\medskip
\begin{flushleft}

 Katsuei  Kenmotsu \\
Mathematical Institute,  Tohoku University  \\
980-8578 \quad  Sendai, Japan \\
email:  kenmotsu@m.tohoku.ac.jp
\end{flushleft}
\end{document}